\documentclass[12pt]{article}

\usepackage[a4paper,hmargin=2cm,vmargin=3cm]{geometry}
\usepackage{fancyhdr,graphicx,lastpage}
\usepackage{listings}
\usepackage{amsmath}
\usepackage{amssymb}
\usepackage{amsthm}
\usepackage{enumerate}
\usepackage{parskip}
\usepackage{multirow}
\usepackage{url}
\usepackage{makeidx}
\usepackage{mathrsfs}
\DeclareMathOperator{\R}{R}
\DeclareMathOperator{\N}{N}
\DeclareMathOperator{\ra}{r}
\DeclareMathOperator{\n}{n}

\newtheoremstyle{krisjan}
{16pt}
{}
{}
{}
{\bf}
{.}
{3pt}
{}

\theoremstyle{krisjan}
\newtheorem{thm}{Theorem}
\newtheorem{prop}[thm]{Proposition}
\newtheorem{lem}[thm]{Lemma}
\newtheorem{cor}[thm]{Corollary}
\title{A note on products of quadratic matrices of singular type.}
\author{C.J. Hattingh\footnote{Corresponding author. krisjan@simulasie.co.za}
}
\begin{document}
\maketitle

\begin{abstract}
I provide an alternate and self-contained proof of Botha's theorem on products of idempotent and square-zero matrices where the product contains two square-zero factors, and provide a conclusive characterisation of products of singular quadratic matrices based on previous results.
\end{abstract}

\section{Introduction}
A singular quadratic matrix which is not the zero matrix has minimum polynomial $x^2-cx$ (where $c$ is some scalar over a field), and is therefore either square-zero, or a scalar multiple of an idempotent matrix. In this article I aim to summarise previous results on products of such matrices in one general statement, presented here as corollary \ref{main_result_conclusion}. 

I want to start however, by presenting an alternate proof to one of the results first proved by Botha \cite{botha_sqi}. In this article I aim to present additional insight into that investigation through a self-contained alternate proof which might provide further insight and aid further research into products of quadratic matrices.

I will now fix some notation. The set of order $n \times n$ square matrices over a field $\mathscr{F}$ is indicated as $M_n(\mathscr{F})$. Matrices are generally indicated by capitals, and vectors in lower case. The null space of a matrix $G$ is indicated as $\N(G)$ and its range as $\R(G)$, the corresponding dimensions of these subspaces are indicated as $\n(G)$ (nullity of $G$) and $\ra(G)$ (rank of $G$) respectively.

{\bf Definition.} $\text{n}_0(G) = \text{n}(G) - \dim(\text{R}(G) \cap \text{N}(G))$.

I denote the simple Jordan block with characteristic value zero and order $k \times k$ as
\[J_k(0) = \begin{bmatrix} \mathbf{0} & 0\\ I_{k-1} & \mathbf{0}\end{bmatrix}.\] It is easy to see that $\n_0(G)$ is the number of blocks $J_1(0)$ in the rational canonical form of $G$.

A vector space is generally indicated in capitals by calligraphic font. The vector space of $n$-dimensional (column) vectors are indicated as $\mathscr{F}^n$,  and the standard basisvector with a 1 in entry $i$ and zeros elsewhere is indicated as $e_i$.

Well known previous results employed in the proof include 
\begin{enumerate}
\item Grassmann's theorem: Let $\mathcal{V}$ be a vector space over a field $\mathscr{F}$ and let $\mathcal{W}$ and $\mathcal{Y}$ be subspaces of $\mathcal{V}$ satisfying the condition that $\mathcal{W}+\mathcal{Y}$ has finite dimension. Then $\dim(\mathcal{W}+\mathcal{Y}) = \dim(\mathcal{W})+\dim(\mathcal{Y})-\dim(\mathcal{W}\cap\mathcal{Y})$.
\item Sylvester's theorem: Let $A,B$ be matrices over $\mathscr{F}$ such that $A$ has $n$ columns and $B$ has $n$ rows. Then (i) $\n(AB) \leq \n(A)+\n(B)$; $\ra(A)+\ra(B)-n \leq \ra(AB) \leq \min\{\ra(A),\ra(B)\}$.
\item Fitting's lemma: $G \in M_n(\mathscr{F})$ is similar to $N \oplus B$ where $N$ is nilpotent and $B$ is non-singular.
\end{enumerate} 

Recent previous results that are employed:
\begin{enumerate}
\item $H \in M_n(\mathscr{F})$ is a product of $k \geq 1$ idempotent matrices (over $\mathscr{F}$) of nullities $n_1, \ldots, n_k$ respectively if and only if $\n(H) \geq n_i \geq 0$ and $n_1 + \cdots + n_k \geq \ra(I-H)$. \cite{botha_idem} \cite{knuppel} 
\item $F \in M_n(\mathscr{F})$ is a product of two square-zero matrices $Z_1, Z_2$ if and only if $\ra(F) \leq \n_0(F)$, and whenever $F$ is such a product the ranks of $Z_1, Z_2$ can be arbitrarily chosen subject to $\ra(G) \leq \ra(Z_1), \ra(Z_2) \leq n/2$. \cite[theorem 3]{botha_sqz}
\end{enumerate}

Finally note that, up to similarity, the ordering of factors in a product consisting entirely of idempotent and square-zero matrices with prescribed nullities may be chosen arbitrarily, due to the fact that the idempotent and square-zero properties are preserved by similarity, and furthermore a product $AB$ is similar to $B^TA^T$ (its transpose) for any two square matrices $A, B$ of the same order. We may therefore generalize the results below to include the case where $G$ consists of any ordering of the factors presented.

\section{Products with two square-zero factors}
\begin{thm} \label{main_result_2}
Let $G \in M_n(\mathscr{F})$ and $n_1, \ldots, n_k, n_{Z_1}, n_{Z_2} \in \mathbb{N}$. The following two statements are equivalent.
\begin{enumerate}
\item $G = E_1\ldots E_k Z_1 Z_2$ where $E_i^2=E_i$ and $\n(E_i) = n_i$ for each $i \in \{1,\ldots, k\}$, and $Z_1^2=Z_2^2=\mathbf{0}$ and $\n(Z_1) = n_{Z_1}$ and $\n(Z_2) = n_{Z_2}$.
\item \begin{enumerate}[(i)] \item $\ra(G) \leq n_1 + \cdots + n_k + \n_0(G)$, \\
 \item $n_1,\ldots, n_k, n_{Z_1}, n_{Z_2} \leq \n(G)$,\\
 \item $n_{Z_1}, n_{Z_2} \geq n/2$. 
 \end{enumerate}
\end{enumerate}
\end{thm}

\begin{lem} \label{rg_nf}
Let $G = HF$ where $H$ and $F$ are arbitrary square matrices. Then $\dim(\text{R}(G) \cap \text{N}(F)) \geq \text{n}(F) - \text{n}_0(G)$.
\end{lem}
{\bf Proof.} \[(\text{R}(G) \cap \text{N}(G)) + \text{N}(F) \subseteq \text{N}(G),\]
and so by Grassmann's theorem 
\begin{equation} \dim(\text{R}(G) \cap \text{N}(G)) + \text{n}(F) - \dim(\text{R}(G) \cap \text{N}(G) \cap \text{N}(F)) \leq \text{n}(G).\label{lem_rg_nf_eq_1} \end{equation} 
Now since $N(F) \subseteq N(G)$, the desired result follows directly from \eqref{lem_rg_nf_eq_1}.\hfill $\square$

\begin{lem} \label{rg_rf_nf}
Let $G = HF$ where $H$ is an arbitrary square matrix and $F=Z_1Z_2$, where $Z_1^2=Z_2^2=\mathbf{0}$. If $\text{R}(G) \cap \text{N}(F) = \text{R}(G) \cap \text{R}(F) \cap \text{N}(F)$ then $\text{r}(G) \leq \text{n}_0(G)$.
\end{lem}
{\bf Proof.} We have 
\[ \text{R}(G) \cap \text{N}(F) = \text{R}(G) \cap \text{R}(F) \cap \text{N}(F) \subseteq \text{R}(F) \cap \text{N}(F),\] and since $F$ is the product of two square-zero matrices we have by \cite[theorem 3]{botha_sqz}
\[\dim(\text{R}(F) \cap \text{N}(F)) \leq \text{n}(F) - \text{r}(F).\] 
It follows by lemma \ref{rg_nf} that
\[\n(F) - \text{n}_0(G) \leq \text{n}(F) - \text{r}(F),\]
and therefore \[\text{r}(G) \leq \text{r}(F) \leq \text{n}_0(G).\]
\hfill $\square$

{\bf Remark:} Notice that the result above includes the case where $\text{R}(G) \cap \text{N}(F) = \{0\}$.
\begin{prop} \label{main_result_2_part_1}
Let $G = HF$ where $H$ is an arbitrary square matrix and $F=Z_1Z_2$, where $Z_1^2=Z_2^2=\mathbf{0}$. Then $\ra(G) \leq \n_0(G) + \ra(I-H)$.\end{prop}
{\bf Proof.} If $\text{R}(G) \cap \text{N}(F) = \text{R}(G) \cap \text{R}(F) \cap \text{N}(F)$ then the result follows directly from lemma \ref{rg_rf_nf}, so suppose that $\text{R}(G) \cap \text{R}(F) \cap \text{N}(F) \subset \text{R}(G) \cap \text{N}(F)$. Then there exists a subspace $\mathcal{W}$ such that $\dim(G(\mathcal{W})) = \dim(\mathcal{W})>0$ and
\begin{eqnarray} 
\nonumber \text{R}(G) \cap \text{N}(F) &=& (\text{R}(G) \cap \text{R}(F) \cap \text{N}(F)) \oplus G(\mathcal{W})\\
\nonumber &\subseteq& (\text{R}(F) \cap \text{N}(F)) + G(\mathcal{W}).
\end{eqnarray}
Now by a similar argument as employed in the proof of lemma \ref{rg_rf_nf} it is easy to show
\[\text{r}(G) \leq \text{n}_0(G) + \dim(G(\mathcal{W})),\] and it remains to show that $\dim(G(\mathcal{W})) \leq \text{r}(I-H)$. 

Let $\{w_1, w_2, \ldots, w_k\}$ be a basis for $\mathcal{W}$, then since $\mathcal{W} \cap \text{N}(G) = \{0\}$ it follows that $\{Gw_1, Gw_2, \ldots, Gw_k\}$ is a linearly independent set, and since $\text{R}(F) \cap G(\mathcal{W}) = \{0\}$ we then also have $\{(F-G)w_1, (F-G)w_2, \ldots, (F-G)w_k\}$ is linearly independent. Now $F-G = (I-H)F$ and it follows that $\{(F-G)w_1, (F-G)w_2, \ldots, (F-G)w_k\} \subseteq \text{R}(I-H)$, which yields the desired result. \hfill $\square$

\begin{lem} \label{jordan_shuffle}
Let $k \geq 2$, then $J_k(0) = EF$ where $\text{N}(F) = \text{N}(J_k(0))$ and $\text{n}_0(F)=1$ and $E^2=E$ and $\text{n}(E) = 1$.
\end{lem}
{\bf Proof.} If $k=2$:
\[\begin{bmatrix} 0 & 0 \\ 1 & 0 \end{bmatrix} = \begin{bmatrix} 0 & 0 \\ 1 & 1 \end{bmatrix} \begin{bmatrix} 1 & 0 \\ 0 & 0 \end{bmatrix}.\]
If $k > 2$, let \[F_1 = \begin{bmatrix} I_{k-2}\\ 0\end{bmatrix},\] then:
\[\begin{bmatrix} \mathbf{0} & 0 \\ I_{k-1} & \mathbf{0} \end{bmatrix} = 
\begin{bmatrix} 0 & \mathbf{0} \\ e_{k-1} & I_{k-1} \end{bmatrix} \begin{bmatrix} e_{k-1}^T & 0 \\ F_1 & \mathbf{0} \end{bmatrix}.\]
\hfill $\square$

{\bf Remark:} Notice that $n_0(J_k(0)) = 0$ for any integer $k\geq 2$.

\begin{cor}
Let $J = J_{k_1}(0) \oplus \cdots \oplus J_{k_m}(0)$ where $k_i \geq 2$. Then $J = EF$ where $E^2=E$ and $\N(F)=\N(J)$ and $0 \leq \text{n}_0(F) \leq \text{n}(E) \leq m$.
\end{cor}
{\bf Proof.}  Let $0 \leq s \leq m$, and let $J_{k_i}(0) = E_iF_i$ as in the lemma above for $1 \leq i \leq s$. Furthermore, choose $E_i = I_{k_i}$ or $E_i=0 \oplus I_{k_i-1}$ for $s < i \leq m$, so that $J_{k_i}(0) = E_i J_{k_i}(0)$ in this case.

Then $J=EF$ with $F=F_1\oplus \cdots \oplus F_s \oplus J_{k_{s+1}}(0) \oplus \cdots \oplus J_{k_{m}}(0)$ and $E = E_1 \oplus \cdots \oplus E_m$. \hfill $\square$

{\bf Proof of theorem \ref{main_result_2}.} Suppose point number 1 in the statement of the theorem is true. Then (ii) and (iii) follow easily, and (i) follows from proposition \ref{main_result_2_part_1} and \cite{botha_idem}.

Suppose point number 2 in the statement of the theorem holds. By the final paragraph of the introduction we may assume without loss of generality that $n_1 \leq \cdots \leq n_k$. 

By Fitting's lemma $G$ is similar to $N \oplus B$ where $N$ is nilpotent and $B$ is invertible. Without loss of generality we can assume that $N = \mathbf{0}_{\text{n}_0(G)} \oplus J$ where $J = J_{k_1}(0) \oplus \cdots \oplus J_{k_m}(0)$ with $k_i \geq 2$ and $m=\dim(\text{R}(G) \cap \text{N}(G))$.

Now $\mathbf{0}_{\text{n}_0(G)} = H_1\mathbf{0}_{\text{n}_0(G)}$ where $H_1$ is any idempotent matrix with $0 \leq \n(H_1) \leq \text{n}_0(G)$. 

Let $s=\max\{\text{r}(G) - \text{n}_0(G),0\}$. It then follows that $s \leq n_1 + \cdots + n_k$ and the corollary above shows that $J=H_2F_2$ where $H_2$ is idempotent, $\N(F_2)=\N(J)$ and $\text{n}_0(F_2)=s \leq \n(H_2) \leq m$.

Finally we have $B = I_{\text{r}(G)}B$.

Now $E = H_1 \oplus H_2 \oplus I_{\text{r}(G)}$ is idempotent. Notice that by the choice of $H_1$ and $H_2$, we can specify $\n(E)$ arbitrarily subject to $s \leq \n(E) \leq \n(G)$ and therefore let $\n(E) = \max\{n_k,s\}$. Now since $\text{n}(E) = \text{r}(I-E)$, the result \cite{botha_idem} shows that $E$ can be written as $E_1\ldots E_k$.

Furthermore $F= \mathbf{0}_{\text{n}_0(G)} \oplus F_2 \oplus B$ is such that $\text{n}_0(F) = \n_0(G)+s \geq \text{r}(G) = \text{r}(F)$. By theorem 3 of \cite{botha_sqz} it follows that $F=Z_1Z_2$.
\hfill $\square$

\section{Products of quadratic matrices of singular type}
\begin{cor} \label{main_result_conclusion}
Let $G \in M_n(\mathscr{F})$ and $n_1, \ldots, n_k, m_{1}, \ldots, m_{l} \in \mathbb{N}$ and $c_1, \ldots, c_k$ be nonzero scalars in $\mathscr{F}$. Set $c = (c_1\ldots c_k)^{-1}$. The following two statements are equivalent.
\begin{enumerate}
\item $G = (c_1E_1)\ldots (c_kE_k) Z_1 \ldots Z_l$ where $E_i^2=E_i$ and $\n(E_i) = n_i$ for each $i \in \{1,\ldots, k\}$, and $Z_j^2=\mathbf{0}$ and $\n(Z_j) = m_{j}$ for each $j  \in \{1,\ldots, l\}$.
\item \begin{enumerate}[(i)] 
\item $n_1, \ldots, n_k, m_{1}, \ldots, m_{l} \leq \n(G)$, 
\item $m_{1}, \ldots, m_{l} \geq n/2$, 
\item furthermore if $l=0$ then $\ra(I-cG) \leq n_1 + \cdots + n_k$,\\ 
if $l=1$ then $\dim(\R(G)+\N(G)) \leq n_1 + \cdots + n_k + m_1$,\\
if $l=2$ then $\ra(G) \leq n_1 + \cdots + n_k + \n_0(G)$.
 \end{enumerate}
\end{enumerate}
\end{cor}
{\bf Proof.} Note that $G = (c_1E_1)\ldots (c_kE_k) Z_1 \ldots Z_l$ is equivalent to $cG = E_1\ldots E_k Z_1 \ldots Z_l$. Now the corollary follows easily by combining the results presented in \cite{botha_idem}, \cite{knuppel}, \cite{botha_sqi}, and the preceding sections, with the fact that for any subspace $S$ we have $cS = \{cv : v \in S\} = S$. \hfill $\square$

\end{document}